\documentclass{elsart}
\usepackage{amssymb}
\usepackage{graphicx}

\begin{document}

\begin{frontmatter}

\title{Geometry~of~Planar~Quadratic~Systems}

\author{Valery A. Gaiko\thanksref{label1}}
\ead{vlrgk@yahoo.com}
\thanks[label1]{The author is grateful to the Department of Applied
Mathematical Analysis of TU~Delft (the Netherlands), the
Department of Mathematical Sciences of the University of Memphis
(USA), and the Institut des Hautes \'{E}tudes Scientifiques
(France) for their hospitality in 2003--2006.}
\address{Department~of~Mathematics
Belarusian~State~University~of~Informatics~and~Radioelectronics
L.\,Beda~Str.\,6-4,~Minsk~220040,~Belarus}

\begin{abstract}
In this paper, the global qualitative analysis of planar quadratic
dynamical systems is established and a new geometric approach to
solving \emph{Hilbert's Sixteenth Problem} in this special case of
polynomial systems is suggested. Using geo\-met\-ric properties of
four field rotation parameters of a new canonical system which is
constructed in this paper, we present a proof of our earlier
conjecture that the maximum number of limit cycles in a~quadratic
system is equal to four and the only possible their distribution
is $(3:1)$~\cite{Gaiko}. Besides, applying the Wintner--Perko
termination principle for multiple limit cycles to our canonical
system, we prove in a different way that a quadratic system has
at~most three limit cycles around a singular point (focus) and
give another proof of the same conjecture.
    \par
    \bigskip
\noindent \emph{Keywords}: Hilbert's sixteenth problem;
Wintner--Perko termination principle; planar quadratic dynamical
system; field rotation parameter; bifurcation; limit cycle;
separatrix cycle

\end{abstract}

\end{frontmatter}

\section{Introduction}
\label{1}

We consider planar dynamical systems
    $$
    \dot{x}=P_{2}(x,y),
    \quad
    \dot{y}=Q_{2}(x,y),
    \eqno(1.1)
    $$
where $P_{2}$ and $Q_{2}$ are quadratic polynomials with real
coefficients in the real variables $x,\:y.$  The main problem of
qualitative theory of such systems is \emph{Hilbert's Sixteenth
Problem} on the maximum number and relative position of their
limit cycles, i.\,e., closed isolated trajectories of (1.1)
\cite{Gaiko}, \cite{ilyash}. In this paper, we suggest a new
geometric approach~\cite{gaib} to studying limit cycle
bifurcations of (1.1) and to solving the \emph{Problem} in this
special case of polynomial systems.
    \par
In particular, in Section~2, we construct two canonical systems
with field rotation parameters, one of which contains four such
parameters. In Section~3, using the canonical systems and
geometric properties of the spirals filling the interior and
exterior domains of limit cycles, we obtain the main result of
this paper on the maximum number and relative position of limit
cycles. In Section~4, we obtain the same result in a different
way: applying the Wintner--Perko termination principle for
multiple limit cycles. Finally, in Conclusion, we discuss some
different approaches to \emph{Hilbert's Sixteenth Problem.}

\section{Canonical systems}
\label{2}

In~\cite{chergai}, \cite{Gaiko}, we constructed a canonical
quadratic system with two field rotation parameters for studying
limit cycle bifurcations:
    $$
    \dot{x}=P(x,y)+\alpha\,Q(x,y),
    \qquad
    \dot{y}=Q(x,y)-\alpha P(x,y),\\[-4mm]
    \eqno(2.1)
    $$
where
    $$
    P(x,y)=-y+mxy+(n-\gamma)y^{2}\!,\quad
    Q(x,y)=x-x^{2}+\gamma\,xy+cy^{2}\!.
    $$
In~\cite{chergai}, \cite{Gaiko}, we show also by which linear
transformations of the phase variables $x,\:y$ arbitrary quadratic
system (1.1) is reduced to form (2.1) and how the parameters of
(2.1) are expressed via the parameters of (1.1). System (2.1) is
especially convenient for the investigation of quadratic systems
in the case two finite singularities when the parameters
$\alpha,\:\gamma$ rotate the field of (2.1) in the whole phase
plane $x,\:y.$
    \par
Later, we constructed a canonical system with three field rotation
parameters, $\alpha,\:\beta,\:\lambda,$
    $$
    \dot{x}=-(1+x)\,y+\alpha\,Q(x,y),
    \quad
    \dot{y}=Q(x,y),
    \eqno(2.2)
    $$
where
    $$
    Q(x,y)=x+\lambda y+ax^{2}+\beta(1+x)y+cy^{2},\\[4mm]
    $$
which, together with the system
    $$
    \dot{x}=-y+\nu y^{2},
    \quad
    \dot{y}=Q(x,y),
    \quad
    \nu=0;1,
    \eqno(2.3)
    $$
can be used in an arbitrary case of finite singularities
\cite{Gaiko}.
    \par
Applying a similar approach, we can construct a canonical system
with the maximum number of field rotation parameters, namely: with
four such para\-me\-ters. It is valid the following theorem.
    \par
    \medskip
\noindent\textbf{Theorem 2.1.}
    \emph{A quadratic system with limit cycles can be reduced to
the canonical form}
    $$
\begin{array}{l}
    \dot{x}=-y\,(1+x+\alpha\,y)\equiv P,\\
    \dot{y}=x+(\lambda+\beta+\gamma)y+a\,x^{2}
        +(\alpha+\beta+\gamma)xy+c\,\gamma\,y^{2}
        \equiv Q
\end{array}
    \eqno(2.4)
    $$
\emph{or}
    \\[-4mm]
    $$
\begin{array}{l}
    \dot{x}=-y\,(1+\nu\,y),
        \quad \nu=0;1,\\
    \dot{y}=x+(\lambda+\beta+\gamma)y+a\,x^{2}
        +(\beta+\gamma)xy+c\,\gamma\,y^{2}.
\end{array}
    \eqno(2.5)
    $$
    \par
\noindent\textbf{Proof.} In \cite{Gaiko} is shown that an
arbitrary quadratic system with limit cycles, by means of Erugun's
two-isocline method~\cite{erug}, can be reduced to the form
    $$
\begin{array}{l}
    \dot{x}=-y+mxy+ny^{2},\\
    \dot{y}=x+\lambda y+ax^{2}+bxy+cy^{2},
\end{array}
    \eqno(2.6)
    $$
where $m=-1$ or $m=0.$
    \par
Input the field rotation parameters into this system so that (2.4)
corresponds to the case of $m=-1$ and (2.5) corresponds to the
case of $m=0.$
    \par
Compare (2.4) with (2.6) when $m=-1.$ Firstly, we have changed
several parameters: $n$ by $-\alpha;$ $b$ by $\beta;$ $c$ by
$c\,\gamma.$ Secondly, we have input additional terms into the
expression for $\dot{y}\!:$ $(\beta+\gamma)\,y$ and
$(\alpha+\gamma)\,xy.$ Similar transformations have been made in
system (2.6) when $m=0;$ but in this case, we have denoted $n$ by
$\nu$ assigning two principal values to this parameter: 0 and 1.
It is obvious that all these transformations do not restrict
generality of systems (2.4) and (2.5) in comparison with system
(2.6), what proves the theorem.
    \medskip
    \par
System (2.4) will be a basic system for studying limit cycle
bifurcations. It contains four field rotation parameters:
$\lambda,\:\alpha,\:\beta,\:\gamma.$ The following lemma is valid
for this system (a similar lemma is valid for system (2.5), with
respect to the parameters $\lambda,\:\beta,\:\gamma).$
    \par
    \medskip
\noindent\textbf{Lemma 2.1.}
    \emph{Each of the parameters $\lambda,\:\beta,\:\gamma,$ and $\alpha$
rotates the vector field of $(2.4)$ in the domains of existence of
its limit cycles, under the fixed other parameters of this system,
namely$\,:$ when the parameter $\lambda,\:\beta,\:\gamma,$ or
$\alpha$ increases $($decreases$),$ the field is rotated in
positive $($negative$)$ direction, i.\,e., counterclockwise
$($clockwise$),$ in the do\-mains, respectively$\,:$}\\[-2mm]
    $$
    1+x+\alpha\,y<0 \ (>0);
    $$
    $$
    (1+x)(1+x+\alpha\,y)<0 \ (>0);
    $$
    $$
    (1+x+c\,y)(1+x+\alpha\,y)<0 \ (>0);
    $$
    $$
    (\lambda+\beta+\gamma)\,y+(a-1)\,x^{2}+
    (\beta+\gamma)\,xy+c\,\gamma\,y^{2}<0 \ (>0).
    $$
    \par
\noindent\textbf{Proof.} \ Using the definition of a field
rotation parameter~\cite{duff},~\cite{Gaiko} we can calculate the
following determinants:\\
    $$
\Delta_{\lambda}=PQ'_{\lambda}-QP'_{\lambda}=-y^{2}(1+x+\alpha\,y);
    $$
    $$
\Delta_{\beta}=PQ'_{\beta}-QP'_{\beta}=-y^{2}(1+x)(1+x+\alpha\,y);
    $$
    $$
\Delta_{\gamma}=PQ'_{\gamma}-QP'_{\gamma}=-y^{2}(1+x+c\,y)(1+x+\alpha\,y);
    $$
    $$
\Delta_{\alpha}=PQ'_{\alpha}-QP'_{\alpha}=~y^{2}((\lambda+\beta+\gamma)\,y+
                (a-1)\,x^{2}+(\beta+\gamma)\,xy+c\,\gamma\,y^{2}).
    $$
    \par
Since, by definition, the vector field is rotated in positive
direction (counter\-clock\-wise) when the determinant is positive
and in negative direction (clock\-wise) when the determinant is
negative~\cite{duff},~\cite{Gaiko} and since the obtained domains
correspond to the domains of existence of limit cycles of (2.4),
the lemma is proved.

\section{The main result}
\label{3}

By means of canonical systems (2.4) and (2.5), we will study
global limit cycle bifurcations of (1.1). First of all, let us
give a new proof of the following theorem.
    \par
    \medskip
\noindent\textbf{Theorem 3.1.}
    \emph{A quadratic system can have at least four limit cycles
in the $(3:1)$-distribution.}
    \medskip
    \par
    \noindent\textbf{Proof.} \
To prove the theorem, consider the case of two finite anti-saddles
and the only saddle at infinity when, for example, $a=1/2$ and
$c=-1$ in (2.4):
    \\[-4mm]
    $$
\begin{array}{l}
    \dot{x}=-y\,(1+x+\alpha\,y),\\
    \dot{y}=x+(\lambda+\beta+\gamma)\,y+(1/2)\,x^{2}
        +(\alpha+\beta+\gamma)\,xy-\gamma\,y^{2}.
\end{array}
    \eqno(3.1)
    $$
Vanish all field rotation parameters in (3.1):
$\alpha=\beta=\gamma=\lambda=0.$ Then we have got a system with
two centers which is symmetric with respect to the $x$-axis.
    \par
Under increasing the parameter $\gamma$ $(0<\gamma\ll1),$ the
vector field of~(3.1) is rotated in negative direction (clockwise)
and the centers turn into foci: $(0,0)$ becomes an unstable focus
and $(-2,0)$ becomes a stable one.
    \par
Fix $\gamma$ and take $\lambda$ satisfying the condition:
$-1\ll\lambda<-\gamma<0$ $(-1\ll\gamma+\lambda<0).$ Then, in the
half-plane $x>-1,$ the vector field of~(3.1) is rotated in
positive direction and the focus $(0,0)$ changes the character of
its stability generating an unstable limit cycle. In the
half-plane $x<-1,$ the field is rotated in negative direction
again and the focus $(-2,0)$ remains stable.
    \par
Fix the parameters $\gamma,$ $\lambda$ and take $\alpha$
satisfying the condition: $\gamma+\lambda\ll\alpha<0.$ After
rotation of the vector field of system~(3.1) in positive
direction, the straight line $x=1$ is destroyed and two limit
cycles are generated by the separatrix cycles formed by this line
and two Poincar\'{e} hemi-circles: a stable limit cycle
surrounding the focus $(0,0)$ and an unstable one surrounding the
focus $(-2,0).$
    \par
Finally, fix the parameters $\gamma,$ $\lambda,$ $\alpha$ and take
$\beta$ satisfying the condition:\linebreak
$0<-\gamma-\lambda<\beta\ll1$
$(0<\beta+\gamma+\lambda\ll-\alpha).$ Then, after rotation of the
vector field in negative direction in the whole phase plane, the
focus $(0,0)$ changes the character of its stability again
generating a stable limit cycle, since the parameter $\alpha$ is
non-rough and negative when $\beta=-\gamma-\lambda.$ Thus, we~have
obtained at least three limit cycles surrounding the focus
$(0,0),$ under the co-existence of a limit cycle surrounding the
focus $(-2,0),$ what proves the theorem.
    \par
It is valid a much stronger theorem (the main result).
    \par
    \medskip
\noindent\textbf{Theorem 3.2.}
    \emph{A quadratic system has at most four limit cycles and
only in the $(3:1)$-distribution.}
    \medskip
    \par
\noindent\textbf{Proof.} \ Consider again the most interesting
case of quadratic systems: with two finite anti-saddles and the
only saddle at infinity when $a=1/2$ and $c=-1$ in (2.4). All
other cases of singular points can be considered in a similar way.
   \par
Vanish all field rotation parameters of system (3.1),
$\alpha=\beta=\gamma=\lambda=0\!:$
    \\[-4mm]
    $$
\begin{array}{l}
    \dot{x}=-y\,(1+x),\\
    \dot{y}=x+(1/2)\,x^{2}.
\end{array}
    \eqno(3.2)
    $$
We have got a system with two centers which is symmetric with
respect to the $x$-axis. Let us input successively the field
rotation parameters into (3.2).
    \par
Begin, for example, with the parameter~$\gamma$ supposing that
$\gamma>0\!:$
    \\[-4mm]
    $$
\begin{array}{l}
    \dot{x}=-y\,(1+x),\\
    \dot{y}=x+\gamma\,y+(1/2)\,x^{2}+\gamma\,xy-\gamma\,y^{2}.
\end{array}
    \eqno(3.3)
    $$
Under increasing $\gamma,$ the vector field of~(3.3) is rotated in
negative direction (clockwise) and the centers turn into foci:
$(0,0)$ becomes an unstable focus and $(-2,0)$ becomes a stable
one.
    \par
Fix $\gamma$ and input a new parameter, for example, $\lambda<0$
into (3.3):
    \\[-4mm]
    $$
\begin{array}{l}
    \dot{x}=-y\,(1+x),\\
    \dot{y}=x+(\lambda+\gamma)\,y+(1/2)\,x^{2}+\gamma\,xy-\gamma\,y^{2}.
\end{array}
    \eqno(3.4)
    $$
Then, in the half-plane $x>-1,$ the vector field of~(3.4) is
rotated in positive direction (counterclockwise) and the focus
$(0,0)$ changes the character of its stability (when
$\lambda=-\gamma)$ generating an unstable limit cycle. Under
decreasing $\lambda,$ this limit cycle will expand until it
disappears in a Poincar\'{e} hemi-cycle with a saddle-node lying
on the invariant straight line $x=-1$~\cite{Gaiko}. In~the
half-plane $x<-1,$ the field is rotated in negative direction
again and the focus $(-2,0)$ remains stable.
    \par
Denote the limit cycle by $\Gamma\!_{1},$ the domain inside the
cycle by $D_{1},$ the domain outside the cycle by $D_{2}$ and
consider logical possibilities of the appearance of other
(semi-stable) limit cycles from a ``trajectory concentration''
surrounding the focus $(0,0).$ It is clear that under decreasing
$\lambda,$ a semi-stable limit cycle cannot appear in the domain
$D_{1},$ since the focus spirals filling this domain will untwist
and the distance between their coils will increase because of the
vector field rotation in positive direction.
    \par
By contradiction, we can also prove that a semi-stable limit cycle
cannot appear in the domain $D_{2}.$ Suppose it appears in this
domain for some values of the parameters $\gamma^{*}>0$ and
$\lambda^{*}<0.$ Return to initial system (3.2) and change the
order of inputting the field rotation parameters. Input first the
parameter $\lambda<0\!:$
    \\[-6mm]
    $$
\begin{array}{l}
    \dot{x}=-y\,(1+x),\\
    \dot{y}=x+\lambda\,y+(1/2)\,x^{2}.
\end{array}
    \eqno(3.5)
    $$
Fix it under $\lambda=\lambda^{*}.$ In the half-plane $x>-1,$ the
vector field of~(3.5) is rotated in negative direction and $(0,0)$
becomes a stable focus. Inputting the parameter $\gamma>0$ into
(3.5), we have got again system (3.4), the vector field of which
is rotated in positive direction in the half-plane $x>-1.$ Under
this rotation, an unstable limit cycle $\Gamma\!_{1}$ will appear
from a Poincar\'{e} hemi-cycle with a saddle-node on the invariant
straight line $x=-1.$ This cycle will contract, the outside
spirals winding onto the cycle will untwist and the distance
between their coils will increase under increasing the parameter
$\gamma$ to the value $\gamma=\gamma^{*}.$ It follows that there
are no values of $\gamma=\gamma^{*}$ and $\lambda=\lambda^{*},$
for which a semi-stable limit cycle could appear in the domain
$D_{2}.$
    \par
This contradiction proves the uniqueness of a limit cycle
surrounding the focus $(0,0)$ in system (3.4) for any values of
the parameters $\gamma$ and $\lambda$ of different signs.
Obviously, if these parameters have the same sign, system (3.4)
has no limit cycles surrounding $(0,0)$ at all, like there are no
limit cycles surrounding the focus $(-2,0)$ for the parameters
$\gamma$ and $\lambda$ of different signs.
    \par
Let system (3.4) have the unique limit cycle $\Gamma\!_{1}.$ Fix
the parameters $\gamma>0,$ $\lambda<0$ and input the third
parameter, $\alpha<0,$ into this system:
    \\[-4mm]
    $$
\begin{array}{l}
    \dot{x}=-y\,(1+x+\alpha\,y),\\
    \dot{y}=x+(\lambda+\gamma)\,y+(1/2)\,x^{2}
    +(\alpha+\gamma)\,xy-\gamma\,y^{2}.
\end{array}
    \eqno(3.6)
    $$
The vector field of~(3.6) is rotated in positive direction again,
the invariant straight line $x=-1$ is immediately destroyed and
two limit cycles appear from the corresponding Poincar\'{e}
hemi-cycles containing this straight line: a stable cycle, denoted
by $\Gamma_{2},$ surrounding the focus $(0,0)$ and an unstable
limit cycle, denoted by $\Gamma_{3},$ surrounding the focus
$(-2,0).$ Under further decreasing $\alpha,$ the limit cycle
$\Gamma_{2}$ will join with $\Gamma\!_{1}$ forming a semi-stable
limit cycle, $\Gamma_{\!12},$ which will disappear in a
``trajectory concentration'' surrounding the origin $(0,0).$ Can
another semi-stable limit cycle appear around the origin in
addition to $\Gamma_{\!12}?$ It is clear that such a limit cycle
cannot appear neither in the domain $D_{1}$ bounded by the origin
and $\Gamma\!_{1}$ nor in the domain $D_{3}$ bounded on the inside
by $\Gamma_{2}$ because of increasing the distance between the
spiral coils filling these domains under decreasing $\alpha.$
    \par
To prove impossibility of the appearance of a semi-stable limit
cycle in the domain $D_{2}$ bounded by the cycles $\Gamma\!_{1}$
and $\Gamma_{2}$ (before their joining), suppose the contrary,
i.\,e., for some set of values of the parameters $\gamma^{*}>0,$
$\lambda^{*}<0,$ and $\alpha^{*}<0,$ such a semi-stable cycle
exists. Return to system (3.2) again and input the parameters
$\alpha<0$ and $\lambda<0\!:$
    \\[-4mm]
    $$
\begin{array}{l}
    \dot{x}=-y\,(1+x+\alpha\,y),\\
    \dot{y}=x+\lambda\,y+(1/2)\,x^{2}+\alpha\,xy.
\end{array}
    \eqno(3.7)
    $$
In the half-plane $x>-1,$ both parameters act in a similar way:
they rotate the vector field of (3.7) in positive direction
turning the origin $(0,0)$ into a stable focus. In the half-plane
$x<-1,$ they rotate the field in opposite directions generating an
unstable limit cycle from the focus $(-2,0).$
    \par
Fix these parameters under $\alpha=\alpha^{*},$
$\lambda=\lambda^{*}$ and input the parameter $\gamma>0$ into
(3.7) getting again system (3.6). Since, on our assumption, this
system has two limit cycles for $\gamma<\gamma^{*},$ there exists
some value of the parameter, $\gamma_{12}$
$(0<\gamma_{12}<\gamma^{*}),$ for which a semi-stable limit cycle,
$\Gamma_{12},$ appears in system (3.6) and then splits into an
unstable cycle, $\Gamma\!_{1},$ and a stable cycle, $\Gamma_{2},$
under further increasing $\gamma.$ The formed domain $D_{2}$
bounded by the limit cycles $\Gamma\!_{1},$ $\Gamma_{2}$ and
filled by the spirals will enlarge, since, on the properties of a
field rotation parameter, the interior unstable limit cycle
$\Gamma\!_{1}$ will contract and the exterior stable limit cycle
$\Gamma_{2}$ will expand under increasing $\gamma.$ The distance
between the spirals of the domain $D_{2}$ will naturally increase,
what will prohibit from the appearance of a semi-stable limit
cycle in this domain for $\gamma>\gamma_{12}.$ Thus, there are no
such values of the parameters, $\gamma^{*}>0,\!$
$\lambda^{*}<0,\!$ $\alpha^{*}<0,$ for which system (3.6) would
have an additional semi-stable limit cycle.
    \par
Obviously, there are no other values of the parameters $\gamma,$
$\lambda,$ $\alpha,$ for which system (3.6) would have more than
two limit cycles surrounding the origin $(0,0)$ and simultaneously
more than one limit cycle surrounding the point $(-2,0)$ (on the
same reasons). It follows that system (3.6) can have at most three
limit cycles and only in the $(2:1)$-distribution.
    \par
Suppose that system (3.6) has two limit cycle, $\Gamma\!_{1}$ and
$\Gamma_{2},$ around the origin $(0,0)$ and the only limit cycle,
$\Gamma\!_{3},$ around the point $(-2,0).$ Fix the parameters
$\gamma>0,$ $\lambda<0,$ $\alpha<0$ and input the fourth
parameter, $\beta>0,$ into (3.6) getting system (3.1). Under
increasing $\beta,$ the vector field of~(3.1) is rotated in
negative direction, the focus $(0,0)$ changes the character of its
stability (when $\beta=-\gamma-\lambda)$ and a stable limit cycle,
$\Gamma_{0},$ appears from the origin. Suppose it happens before
the cycle $\Gamma\!_{1}$ disappears in $(0,0)$ (this is possible
by Theorem~3.1). Under further increasing $\beta,$ the cycle
$\Gamma_{0}$ will join with $\Gamma\!_{1}$ forming a semi-stable
limit cycle, $\Gamma_{\!01},$ which will disappear in a
``trajectory concentration'' surrounding the origin $(0,0);$ the
other cycles, $\Gamma_{2}$ and $\Gamma_{3},$ will expand tending
to Poincar\'{e} hemi-cycles with the straight line $x=-1.$
    \par
Let system (3.1) have four limit cycles: $\Gamma_{0},$
$\Gamma\!_{1},$ $\Gamma_{2}$, and $\Gamma_{3}.$ Can an additional
semi-stable limit cycle appear around the origin under increasing
the parameter $\beta\,?$ It is clear that such a limit cycle
cannot appear neither in the domain $D_{0}$ bounded by the origin
and $\Gamma_{0}$ nor in the domain $D_{2}$ bounded by
$\Gamma\!_{1}$ and $\Gamma_{2}$ because of increasing the distance
between the spiral coils filling these domains under increasing
$\beta.$ Consider two other domains: $D_{1}$ bounded by the cycles
$\Gamma_{0},$ $\Gamma\!_{1}$ and $D_{3}$ bounded on the inside by
the cycle $\Gamma_{2}.$ As before, we will prove impossibility of
the appearance of a semi-stable limit cycle in these domains by
contradiction.
    \par
Suppose that for some set of values of the parameters,
$\gamma^{*}>0,$ $\lambda^{*}<0,$ $\alpha^{*}<0,$ and
$\beta^{*}>0,$ such a semi-stable cycle exists. Return to system
(3.2) again and input first the parameters $\beta>0,$ $\gamma>0$
and then the parameter $\alpha<0\!:$
    \\[-4mm]
    $$
\begin{array}{l}
    \dot{x}=-y\,(1+x+\alpha\,y),\\
    \dot{y}=x+(\beta+\gamma)\,y+(1/2)\,x^{2}
        +(\alpha+\beta+\gamma)\,xy-\gamma\,y^{2}.
\end{array}
    \eqno(3.8)
    $$
Fix the parameters $\beta,$ $\gamma$ under the values $\beta^{*},$
$\gamma^{*},$ respectively. Under decreasing the parameter
$\alpha,$ two limit cycles immediately appear from Poincar\'{e}
hemi-cycles with the straight line $x=-1\!\!:$ a stable cycle,
$\Gamma_{2},$ around $(0,0)$ and an unstable one, $\Gamma_{3},$
around $(-2,0).$ Fix $\alpha$ under the value $\alpha^{*}$ and
input the parameter $\lambda<0$ into (3.8) getting system (3.1).
    \par
Since, on our assumption, system (3.1) has three limit cycles
around the origin $(0,0)$ for $\lambda>\lambda^{*},$ there exists
some value of the parameter, $\lambda_{01}$
$(\lambda^{*}<\lambda_{01}<0),$ for which a semi-stable limit
cycle, $\Gamma_{01},$ appears in this system and then splits into
a stable cycle, $\Gamma_{0},$ and an unstable cycle,
$\Gamma\!_{1},$ under further decreasing $\lambda.$ The formed
domain $D_{1}$ bounded by the limit cycles $\Gamma_{0},$
$\Gamma\!_{1}$ and also the domain $D_{3}$ bounded on the inside
by the limit cycle $\Gamma_{2}$ will enlarge and the spirals
filling these domains will untwist excluding a possibility of the
appearance of a semi-stable limit cycle there, i.\,e., at most
three limit cycles can exist around the origin $(0,0).$ On the
same reasons, a semi-stable limit cannot appear around the point
$(-2,0)$ under decreasing the parameter $\lambda,$ i.\,e., at most
one limit cycle can exist around this point simultaneously with
three limit cycles surrounding $(0,0).$
    \par
All other combinations of the parameters $\lambda,$ $\alpha,$
$\beta,$ $\gamma$ are considered in a similar way. It follows that
system (3.1) has at most four limit cycles and only in the
$(3:1)$-distribution. Applying the same approach to canonical
system (2.5), we can complete the proof of the theorem.

\section{The Wintner--Perko termination principle}
\label{4}

For the global analysis of limit cycle bifurcations in
\cite{Gaiko}, we used the Wintner--Perko termination principle
which was stated for relatively prime, pla\-nar, analytic systems
and which connected the main bifurcations of limit
cycles~\cite{Perko}, \cite{wint}. Let us formulate this principle
for the polynomial system
    $$
    \mbox{\boldmath$\dot{x}$}=\mbox{\boldmath$f$}
    (\mbox{\boldmath$x$},\mbox{\boldmath$\mu$)},
    \eqno(4.1_{\mbox{\boldmath$\mu$}})
    $$
where $\mbox{\boldmath$x$}\in\textbf{R}^2;$ \
$\mbox{\boldmath$\mu$}\in\textbf{R}^n;$ \
$\mbox{\boldmath$f$}\in\textbf{R}^2$ \ $(\mbox{\boldmath$f$}$ is a
polynomial vector function).
    \par
    \medskip
\noindent\textbf{Theorem 4.1
    (Wintner--Perko termination principle).}
    \emph{Any one-para\-me\-ter fa\-mi\-ly of multiplicity-$m$
limit cycles of relatively prime polynomial system \
$(4.1_{\mbox{\boldmath$\mu$}})$ can be extended in a unique way to
a maximal one-parameter family of multiplicity-$m$ limit cycles of
\ $(4.1_{\mbox{\boldmath$\mu$}})$ which is either open or cyclic.}
    \par
\emph{If it is open, then it terminates either as the parameter or
the limit cycles become unbounded; or, the family terminates
either at a singular point of \ $(4.1_{\mbox{\boldmath$\mu$}}),$
which is typically a fine focus of multiplicity~$m,$ or on a
$($compound$\,)$ separatrix cycle of \
$(4.1_{\mbox{\boldmath$\mu$}}),$ which is also typically of
multiplicity~$m.$}
    \medskip
    \par
The proof of the Wintner--Perko termination principle for general
polynomial system $(4.1_{\mbox{\boldmath$\mu$}})$ with a vector
parameter $\mbox{\boldmath$\mu$}\in\textbf{R}^n$ parallels the
proof of the pla\-nar termination principle for the system
    $$
    \vspace{1mm}
    \dot{x}=P(x,y,\lambda),
        \quad
    \dot{y}=Q(x,y,\lambda)
    \eqno(4.1_{\lambda})
    \vspace{2mm}
    $$
with a single parameter $\lambda\in\textbf{R}$ (see \cite{Gaiko},
\cite{Perko}), since there is no loss of generality in assuming
that system $(4.1_{\mbox{\boldmath$\mu$}})$ is parameterized by a
single parameter $\lambda;$ i.\,e., we can assume that there
exists an analytic mapping $\mbox{\boldmath$\mu$}(\lambda)$ of
$\textbf{R}$ into $\textbf{R}^n$ such that
$(4.1_{\mbox{\boldmath$\mu$}})$ can be written as
$(4.1\,_{\mbox{\boldmath$\mu$}(\lambda)})$ or even
$(4.1_{\lambda})$ and then we can repeat everything, what had been
done for system $(4.1_{\lambda})$ in~\cite{Perko}. In particular,
if $\lambda$ is a field rotation parameter of $(4.1_{\lambda}),$
it is valid the following Perko's theorem on monotonic families of
limit cycles.
    \par
    \medskip
\noindent\textbf{Theorem 4.2.}
    \emph{If $L_{0}$ is a nonsingular multiple limit cycle of
$(4.1_{0}),$ then  $L_{0}$ belongs to a one-parameter family of
limit cycles of $(4.1_{\lambda});$ furthermore\/$:$}
    \par
1)~\emph{if the multiplicity of $L_{0}$ is odd, then the family
either expands or contracts mo\-no\-to\-ni\-cal\-ly as $\lambda$
increases through $\lambda_{0};$}
    \par
2)~\emph{if the multiplicity of $L_{0}$ is even, then $L_{0}$
be\-fur\-cates into a stable and an unstable limit cycle as
$\lambda$ varies from $\lambda_{0}$ in one sense and $L_{0}$
dis\-ap\-pears as $\lambda$ varies from $\lambda_{0}$ in the
opposite sense; i.\,e., there is a fold bifurcation at
$\lambda_{0}.$}
    \par
    \medskip
Using Theorems~4.1 and~4.2 in \cite{gai1}--\cite{gai4}, we proved
a theorem on three limit cycles around a singular point for
canonical systems (2.2) and (2.3). Let us prove the same theorem
using systems (2.4) and (2.5).
    \par
    \medskip
\noindent\textbf{Theorem 4.3.}
    \emph{There exists no quadratic system having a swallow-tail
bifurcation surface of multiplicity-four limit cycles in its
pa\-ra\-meter space. In other words, a quadratic system cannot
have neither a multi\-plicity-four limit cycle nor four limit
cycles around a singular point $($focus$\,),$ and the maximum
multi\-plicity or the maximum number of limit cycles surrounding a
focus is equal to three.}
    \medskip
    \par
\noindent\textbf{Proof.} \ The proof of this theorem is carried
out by contradiction. Consider canonical systems (2.4) and (2.5),
where system (2.5) represents two limit cases of (2.4).
    \par
Suppose that system (2.4) with four field rotation parameters,
$\lambda,$ $\alpha,$ $\beta,$ and $\gamma,$ has four limit cycles
around the origin (system (2.5) is considered in a similar way).
Then we get into some domain of the field rotation parameters
being restricted by definite con\-di\-tions on two other
parameters, $a$ and $c,$ corresponding to one of the six cases of
finite singularities which we considered in \cite{Gaiko}. Without
loss of generality, we can fix both of these parameters. Thus,
there is a domain bounded by three fold bifurcation surfaces
forming a swallow-tail bifurcation surface of multiplicity-four
limit cycles in the space of the field rotation pa\-ra\-me\-ters
$\lambda,$ $\alpha,$ $\beta,$ and $\gamma.$
    \par
The cor\-res\-pon\-ding maximal one-parameter family of
multiplicity-four limit cycles cannot be cyclic, otherwise there
will be at least one point cor\-res\-pon\-ding to the limit cycle
of multi\-pli\-ci\-ty five (or even higher) in the parameter
space. Extending the bifurcation curve of multi\-pli\-ci\-ty-five
limit cycles through this point and parameterizing the
corresponding maximal one-parameter family of
multi\-pli\-ci\-ty-five limit cycles by a field-rotation
para\-me\-ter, according to Theorem~4.2, we will obtain a
monotonic curve which, by the Wintner--Perko termination principle
(Theorem~4.1), terminates either at the origin or on some
separatrix cycle surrounding the origin. Since we know absolutely
precisely at least the cyclicity of the singular point (Bautin's
result \cite{BL}) which is equal to three, we have got a
contradiction with the termination principle stating that the
multiplicity of limit cycles cannot be higher than the
multi\-pli\-ci\-ty (cyclicity) of the singular point in which they
terminate.
    \par
If the maximal one-parameter family of multiplicity-four limit
cycles is not cyclic, on the same principle (Theorem~4.2), this
again contradicts to Bautin's result not admitting the
multiplicity of limit cycles higher than three. This contradiction
completes the proof.
    \par
    \medskip
As was shown in \cite{Gaiko}, to complete the solution of
\emph{Hilbert's Sixteenth Problem} for quadratic systems (1.1), it
is sufficient to prove impossibility of the $(2:2)$-distribution
of limit cycles only in the case of two finite foci and a saddle
at infinity. In \cite{chergai} (see also~\cite{Gaiko}), using
canonical system (2.1) with two field rotation parameters,
$\alpha$ and $\gamma,$ in the case of two foci and a saddle at
infinity, we constructed a quadratic system with at least four
limit cycles in the $(3:1)$-distribution. If to let this system
have only three limit cycles in the $(2:1)$-distribution, i.\,e.,
two cycles around the focus $(0,0)$ and the only one around the
focus $(1,0),$ it is easy to show impossibility of obtaining the
second limit cycle around $(1,0)$ by means of the parameters
$\alpha$ and $\gamma.$ Logically, we can suppose only that a
semi-stable cycle appears around the focus $(1,0)$ under the
variation of a field rotation parameter, for example, $\alpha.$
Then, applying the Wintner--Perko termination principle
(Theorem~4.1), we can show that the maximal one-parameter family
of multiplicity-three limit cycles parameterized by another field
rotation parameter, $\gamma,$ cannot terminate in the focus
$(1,0),$ since it will be a rough focus for any $\alpha\neq0$ (see
\cite{chergai}, \cite{Gaiko}). The same proof could be given for
canonical system~(2.4). Thus, we have given one more proof of
Theorem~3.2 on at most four limit cycles in the only
$(3:1)$-distribution for quadratic systems (1.1).

\section{Conclusion}
\label{5}

In~\cite{Gaiko}, applying the methods of catastrophe theory and
the Wintner--Perko termination principle for multiple limit
cycles, we have developed the global bifurcation theory of planar
polynomial dynamical systems and, basing on this theory, we have
suggested a program on the complete solution of \emph{Hilbert's
Sixteenth Problem} for the case of quadratic systems. In
principal, the program has been realized in \cite{Gaiko} (see the
previous section). In this paper, we have presented a new
(geometric) approach to its realization.
    \par
Our program is an alternative to the program which is put forward
in~\cite{dumrr1}, \cite{dumrr2} and which is often called as
``\,Roussarie's program'' by the name of its ideological
inspirer~\cite{Rous}. Roussarie's program is reduced to the
classification of se\-pa\-rat\-rix cycles, determining their
cyclicity and finding an upper bound of the number of limit cycles
for quadratic systems. Unfortunately, there are some serious
problems in the realization of this program: for example, it is
not clear how to determine the cyclicity of non-monodromic
separatrix cycles when there is no return map in the neighborhood
of these cycles and there is no general approach to the study of
the cyclicity of separatrix cycles in the case of center when the
return map is identical zero. Besides, even in the case of
realization of the program, as its authors note
themselves~\cite{dumrr1}, the obtained upper bound of the number
of limit cycles obviously can not be optimal, since the used pure
analytic me\-thods cannot ensure neither the global control of
limit cycle bifurcatins around a singular point nor, especially,
the simultaneous control of the bifurcations around different
singular points.
    \par
Thereupon, it makes sense to say some words on Roussarie's review
MR2023976 (2005d:37102) on \cite{Gaiko}. The only concrete remark
in this ``awkward'' review is the following: ``I just mention the
hazardous claim made in Theorem 4.12, page 137, that there exists
no quadratic system having a swallow-tail bifurcation surface of
multiplicity-four limit cycles. Looking at the proof, it seems
that the author unfortunately confuses two different notions:
paths of limit cycles, as defined in Definition~4.7, page~112, and
lines of multiple limit cycles, as defined by Perko (and recalled
in Definition~4.13, page~127). In fact, there is nothing
forbidding that a path begin at a parameter value with a
multiplicity-four limit cycle and end at a focus point''. So,
Roussarie's remark is related to a swallow-tail bifurcation
surface of multi\-pli\-city-four limit cycles. However, Definition
4.13, page~127, is a definition of a cusp bifurcation surface of
multiplicity-three limit cycles. This is an evident lack of
correspondence! Maybe, the reviewer means Definition~4.14,
pages~128-129? Then it seems that he did not pay attention for our
remark on page~132, following just after Theorem~4.10, which could
perhaps settle his doubts. Moreover, there is a reference to the
corresponding work by Perko in this remark (see also
\cite{Perko}). Or the reviewer has complaints against Perko's
work, too? Besides, his ``claim'' that ``there is nothing
forbidding that a path begin at a parameter value with a
multiplicity-four limit cycle and end at a focus point'' says that
he unfortunately does not see (or does not want to see) Bautin's
result \cite{BL} (Theorem~2.1, page~45) on the cyclicity of a
singular point of the focus or center type, which is an obstacle
on such a path. Or, maybe, Bautin's result is also
``questionable''?

\end{document}